\documentclass[11pt]{amsart}

\usepackage{amssymb,amsmath,amsthm,amsfonts, xypic}
\usepackage[colorlinks, linktocpage=true]{hyperref}
\hypersetup{citecolor=blue}
\usepackage{todonotes}
\usepackage{booktabs, comment}
\usepackage[section]{placeins}

\newtheorem{theorem}{Theorem}[section]

\newtheorem{lemma}[theorem]{Lemma}

\newtheorem{question}[theorem]{Question}

\newcommand{\C}{\mathbb C}
\newcommand{\pc}{{\mathbb P^1(\C)}}
\newcommand{\V}{\mathcal V}

\newcommand{\R}{\mathbb R}
\newcommand{\Q}{\mathbb Q}
\newcommand{\Z}{\mathbb Z}
\newcommand{\h}{\mathbb H}
\def\slc{{\rm SL_2(\C)}}

\def\pslc{{\rm PSL_2(\C)}}
\def\hom{{\rm Hom}}
\def\tr{{\rm tr}}

\title{An Introduction to A-polynomials and their Mahler Measures}

\author{Mehmet Haluk \c{S}eng\"{u}n}
\email{M.H.Sengun@warwick.ac.uk}
\urladdr{http://warwick.ac.uk/haluksengun}
\address{Mathematics Institute, University of Warwick, Coventry, UK}
\thanks{The author is supported by a Marie Curie Intra-European Fellowship.}

\begin{document}

\maketitle
\tableofcontents
 
\section{Introduction}
These are the notes of the three lectures I delivered at the mini-workshop ``Knot Theory and Number Theory around the A-Polynomial" 
at the Instituto Superior T\'ecnico (IST) in Lisbon in January 2014. The goal of the lectures was to familiarize, both the author and, the audience with 
the $A$-polynomials and the connection between the Mahler measures of A-polynomials and volumes. The style of these notes is expository, written 
informally with the aim of giving a flavor of the subject with ample amount of references to direct the interest readers to the details. 

I would like to thank Nuno Freitas for initiating this meeting and to Roger Picken for hosting me. Finally, I would like to thank the audience and the IST for 
providing a warm and welcoming atmosphere. Finally, I thank Wadim Zudilin for helpful correspondence. 

\section{Algebraic sets related to $3$-manifolds}

\subsection{Representation Variety} Let $\Gamma$ be a finitely presented group 
$$\Gamma = \langle g_1, \hdots ,g_n \mid r_1(g_1, \hdots, g_n) = \hdots =r_k(g_1, \hdots, g_n) = e \rangle.$$
The {\bf representation variety} of $\Gamma$, denoted $R(\Gamma)$,  is the set of homomorphisms from $\Gamma$ into $\slc$;
$$R(\Gamma) := \hom( \Gamma, \slc).$$
The map 
$$\rho \mapsto (\rho(g_1), \hdots , \rho(g_n))$$ 
is an injection from $R(\Gamma)$ to $\slc^n$ and its image is in the solution set of equations 
\begin{equation} \label{eqn: matrix-eqns}
R_j(x_1, \hdots , x_n)= 1, \ \ \ \ 1 \leq j \leq k \end{equation} 
in $\slc^n$ coming from the relations $r_1, \hdots, r_k$ of $\Gamma$. 
Conversely, any element in $\slc^n$ solving the equations in (\ref{eqn: matrix-eqns}) corresponds to an element of $R(\Gamma)$. 
Let us put
$$ \rho(g_i) = \begin{pmatrix} a_i & b_i \\ c_i & d_i \end{pmatrix}, \ \ \ 1 \leq i \leq n.$$
As the determinant of each $\rho(g_i)$ is one, we get $n$ equations in the $4n$ variables $a_i,b_i,c_i,d_i$ with $1 \leq i \leq n$. 
Moreover, each matrix equation in (\ref{eqn: matrix-eqns}) gives 4 polynomial equations (all defined over $\Z$) in the same $4n$ variables after replacing every occurrence of 
$\left ( \begin{smallmatrix} a_i & b_i \\ c_i & d_i \end{smallmatrix} \right ) ^{-1}$ with $\left ( \begin{smallmatrix} d_i & -b_i \\ -c_i & a_i \end{smallmatrix} \right )$. 
Thus $R(\Gamma)$ can be identified with an {\bf algebraic set} in $\C^{4n}$, which in general is not irreducible, given by $4k+n$ equations.  

If we work with another presentation for $\Gamma$, then we end up with changing $R(\Gamma)$ by an isomorphism of algebraic sets. Indeed, writing the new generators as words in the old generators translates to making (polynomial) change of variables at the level of algebraic sets. Thus $R(\Gamma)$ is unique up to isomorphism. 

\subsection{Character Variety} 

Recall that the {\em character} of a representation $\rho \in R(\Gamma)$ is the homomorphism $\chi_\rho : \Gamma \rightarrow \C$ defined by $\chi_\rho(g) = \tr(\rho(g))$ for $g \in \Gamma$. The {\bf character variety} of $\Gamma$, denoted $X(\Gamma)$, is the space of characters of elements in $R(\Gamma)$, that is, 
$$X(\Gamma) := \{ \chi : \Gamma \rightarrow \C \mid \chi = \chi_\rho \ \  \textrm{for some} \ \rho \in R(\Gamma) \}.$$
It follows, with more effort than in the case of $R(\Gamma)$, that $X(\Gamma)$ is an affine algebraic set as well.  

Another way to look at $X(\Gamma)$ is to consider the conjugation action of $\slc$ on $R(\Gamma)$. It can be shown that if $\rho \in R(\Gamma)$ is irreducible, then $\rho,\rho'$ are conjugate if and only if $\chi_\rho=\chi_{\rho'}$.  However non-conjugate reducible representations can have the same character.  Thus $X(\Gamma)$ is not 
quite the set-theoretic quotient $R(\Gamma)/\slc$ but categorical quotient $R(\Gamma)// \slc$. 

It can be shown, with a proof based on the identity $\tr(AB)+\tr(AB^{-1})=\tr(A)\tr(B)$, that an element $\chi$ of $X(\Gamma)$ is determined by its values on the elements of $\Gamma$ of the form $g_{i_1}\cdots g_{i_m}$ with $1 \leq m \leq n$ and $1 \leq i_1 < \hdots < i_m \leq n$. Note that there are $2^n-1$ such elements. In fact it can be 
shown that it suffices to consider a set whose size is only polynomial in $n$, see \cite{florentino}. Note that this number also gives an upper bound on the dimension of any component of $X(\Gamma)$. For future use, let us remark that if $\Gamma$ is generated by two elements $a,b$, $\chi(a),\chi(b),\chi(ab^{\pm1})$ suffice to uniquely determines $\chi \in X(\Gamma)$. 

We call a $\chi_\rho \in X(\Gamma)$ {\bf reducible} if $\rho$ is reducible, that is, all the elements in the image $\rho(\Gamma)$ share a common eigenvector. This is equivalent to saying that, up to conjugation, $\rho(\Gamma)$ lies in $ ( \begin{smallmatrix} \star & \star \\ 0 & \star \end{smallmatrix} )$. Otherwise, we call $\chi_\rho$ {\bf irreducible}. 
It can be seen that the subset of reducible characters $X_{{\rm red}}(\Gamma)$ is a sub-algebraic variety of $X(\Gamma)$. We define $X_{{\rm irr}}(\Gamma)$ as the Zariski 
closure of the complement of $X_{{\rm red}}(\Gamma)$ in $X(\Gamma)$. 

Let $M$ be a hyperbolic 3-manifold of finite volume with holonomy representation $\rho_0 : \pi_1M \rightarrow \pslc$. Mostow Rigidity says that $\rho_0$ is unique 
up to conjugation. An irreducible component of $X(M):=X(\pi_1M)$ is called {\bf canonical}, and denoted $X_0(M)$, if it contains the character of a 
lift\footnote{It is known that lifts exist and are parametrized by $H^1(M,\Z / 2\Z) \simeq \Z / 2\Z$.} 
$\rho: \pi_M \rightarrow \slc$ of $\rho_0$. A theorem of Thurston implies that 
$$\dim_\C X_0(M) = \# \{ \textrm{cusps of} \ M \} . $$

\subsection{Examples}
\begin{itemize} 
\item ({\bf 2-torus}) Let $\Gamma \simeq \Z \times \Z \simeq \langle a ,b \mid [a,b]=e \rangle$. Since $\Gamma$ is abelian, for any $\rho \in R(\Gamma)$, 
its image $\rho$ will be so too. Any abelian subgroup of $\slc$ can be conjugated so that it is upper-triangular. So let us put 
$$ \rho(a) = \begin{pmatrix} x & \star \\ 0 & x^{-1} \end{pmatrix}, \ \ \ \ \textrm{and} \ \ \ \ 
\rho(b) = \begin{pmatrix} y & \star \\ 0 & y^{-1} \end{pmatrix}.$$
Clearly $\chi_\rho(a)=x+x^{-1}$, $\chi_\rho(b)=y+y^{-1}$ and $\chi_\rho(ab)=xy+(xy)^{-1}$. Thus we can parametrize $\chi_\rho$ with 
the tuple $(x,y)$. Noting that the pair 
$(x^{-1}, y^{-1})$ gives rise to the same character,  we can identify $X(\Gamma)$ with $\C^\star \times \C^\star$ modulo the involution $(x,y) \mapsto (x^{-1}, y^{-1})$.
\vspace{.1 in}
\item ({\bf Figure 8-knot}) Let $N$ be the compact 3-manifold given by the complement inside the 3-sphere of an open tubular neighborhood of the figure 8-knot. 
It is well-known that the fundamental group $\Gamma$ of $N$ admits the following presentation 
$$\Gamma \simeq \langle a ,b \mid wa=bw, \ \ w=a^{-1}bab^{-1} \rangle.$$
We see from the presentation that for any $\rho \in R(\Gamma)$, $\rho(a)$ and $\rho(b)$ are conjugate and thus 
\begin{equation} \label{equaltrace}
\chi_\rho(a)=\chi_\rho(b). 
\end{equation}
To determine $\chi_\rho$, it suffices to determine $\chi_\rho(ab^{-1}).$ 

Let us first assume that $\rho$ is reducible. In light of Equation \ref{equaltrace}, we put
$$ \rho(a) = \begin{pmatrix} m & \star \\ 0 & m^{-1} \end{pmatrix} \ \ \ \ \textrm{and} \ \ \ \ 
\rho(b) = \begin{pmatrix} m & \star \\ 0 & m^{-1} \end{pmatrix}$$
and see that 
$$\rho(ab)=\begin{pmatrix} 1 & \star \\ 0 & 1 \end{pmatrix}.$$
Thus $\chi_\rho(ab^{-1})=2$ for any reducible $\rho \in \R(\Gamma)$ implying that $\chi_\rho$ for which $\rho$ is reducible can be simply identified with $\{ (m+m^{-1}, 2) \} \simeq \C$.  

Let us now assume that $\rho$ is irreducible. We will need the following lemma from \cite[Lemma 7]{riley}. 
\begin{lemma} Let $M_1,M_2$ be noncommuting elements of $\slc$ with the same trace. 
Then there exist $t,u$ such that
$$U M_1 U^{-1} =  \begin{pmatrix} t & 1 \\ 0 & t^{-1} \end{pmatrix}, \ \ 
U M_2 U^{-1} = \begin{pmatrix} t & 0 \\ -u & t^{-1} \end{pmatrix},$$
where $U=\left ( \begin{smallmatrix} u & 0 \\ 0 & u^{-1} \end{smallmatrix} \right )$. If $F$ is the field generated by the entries of $M_1$ and $M_2$, then 
$t,u$ belong to an at most quadratic extension of $F$.
\end{lemma}
\end{itemize}
Following the above lemma and \cite{petersen}, we put
$$ \rho(a) = \begin{pmatrix} t & 1 \\ 0 & t^{-1} \end{pmatrix}, \ \ \ \ \textrm{and} \ \ \ \ 
\rho(b) = \begin{pmatrix} t & 0 \\ 2-u & t^{-1} \end{pmatrix}.$$
We see that $\rho(ab^{-1})$ has trace\footnote{The choice of $``2-u"$ as opposed to $``-u"$ as in the Lemma is inessential and 
is only made to make the trace look nicer.} $u$. So we see that $\chi_\rho$ is parametrized by two parameters $x:=t+t^{-1}$ and $u$. 
The question now is whether there are any relations between these parameters. To answer this question, we use the equation 
$$\rho(w)\rho(a)-\rho(b)\rho(w)= \begin{pmatrix} 0 & 0 \\ 0 & 0 \end{pmatrix}$$
coming from the relation 
$``wa=bw"$. Direct computation show that we have $\rho(w)\rho(a)-\rho(b)\rho(w)$ is equal to
$$\begin{pmatrix} 0 & X  \\  
(u-2)X & 0 \end{pmatrix}$$
where $Z=(x^2-2)(1-u)+1-u+u^2$. Thus the $\chi_\rho$ for which $\rho$ is irreducible is in bijection with the algebraic set 
$$\{ (x,u) \in \C^2 \mid (x^2-2)(1-u)+1-u+u^2 = 0 \}.$$
Manipulating this as 
$$(x^2-2)(1-u)+1-u+u^2 = 0  \Leftrightarrow  x^2(u-1)=u^2+u-1 \Leftrightarrow  (x(u-1))^2 =(u^2+u-1)(u-1) $$
and making the change of variable $z=x(u-1)$, we arrive at the curve
$$z^2 = u^3-2u+1$$
which is the elliptic curve of Cremona label ``40a3". 
\subsection{The A-polynomial} Let $N$ be a compact 3-manifold with boundary a torus $T$. 
In this section we will talk about the A-polynomial of $N$. The idea is that instead of working with $X(\pi_1 N)$, which can be quite big, 
we work with a certain plane curve $D_N$, which still carries a great deal information about $N$. Roughly speaking the curve $D_N$ arises 
as the image of $X(\pi_1N)$ in $X(\pi_1T)$ which arises from restriction of representation of $\pi_1N$ to $\pi_1T \leq \pi_1N$. The A-polynomial 
generates the defining ideal of the curve $D_N$.

Let us give details. Put $\Gamma= \pi_1N$. Consider the subset $R_T(\Gamma)$ of $R(\Gamma)$ formed by those $\rho$ for which $\rho(\pi_1T)$ is 
upper-triangular. Note that this subset is algebraic as one can polynomially express being upper-triangular. 
Let us fix a basis $\pi_1T \simeq \Z \times \Z = <m,\ell>$. Putting 
$$ \rho(m) = \begin{pmatrix} M & \star \\ 0 & 1/M \end{pmatrix}, \ \ \ \ \textrm{and} \ \ \ \ 
\rho(\ell) = \begin{pmatrix} L & \star \\ 0 & 1/L \end{pmatrix},$$
we form the {\em eigenvalue map}
$$\varepsilon : R_T(\Gamma) \rightarrow \C^2$$
defined by $\rho \mapsto (M,L)$. Let $\V$ denote the Zariski closure of the image of $\varepsilon$ in $\C^2$. It turns out that 
all the components of $\V$ are either zero or one dimensional, see \cite[Lemma 2.1]{DG}. Each one dimensional component is the zero set of a single polynomial (as they are hypersurfaces) with two variables and so for each such component $\mathcal{C}_i$ of $\V$, we fix such a polynomial $c_i(M,L)$ and define the {\bf A-polynomial} of $N$ 
to be 
$$A_N(M,L) := \prod_{\mathcal{C}_i} c_i(M,L)$$ 
where we run over all the one-dimensional components $\mathcal{C}_i$ of $\V$. 
Note that $\V$ differs than $\varepsilon(R_T(\Gamma))$ by at most finitely many points. So with at most finitely many exceptions, a point $(x_0,y_0)$ satisfies $A_N(M,L)$ 
if and only if there is $\rho \in R(\Gamma)$ such that 
$$ \rho(m) = \begin{pmatrix} x_0 & \star \\ 0 & x_0^{-1} \end{pmatrix}, \ \ \ \ \textrm{and} \ \ \ \ 
\rho(\ell) = \begin{pmatrix} y_0 & \star \\ 0 & y_0^{-1} \end{pmatrix}.$$
The A-polynomial depends on our choice of basis for $\pi_1 \partial N$. Changing the basis results in multiplying the A-polynomial by powers of $M$ and $L$. 
Moreover, it is easy to observe that $A_N(M,L)=A_N(1/M,1/L)$.

Another way to approach the A-polynomial is through the restriction map. The inclusion $i : T \hookrightarrow N$ induces a regular map 
of complex algebraic sets $i^\star : X(N) \rightarrow X(T)$ via considering the restrictions of the representations of $\Gamma$ to $\pi_1T$. 
Let $V$ be the 1-dimensional part of the image $i^\star(X(N))$. Then A-polynomial is the defining polynomial of the plane curve $\V$ in $\C^* \times \C^*$ 
obtained by lifting $V$ under the surjection $\C^* \times \C^* \rightarrow X(T)$. Since all the maps involved are over $\Q$, we can arrange things so that the 
A-polynomial has integer coefficients.

The naive way of computing the A-polynomial goes as follows, see \cite{CL-remarks}. Assume that we have identified $R_T(\Gamma)$ as an algebraic subset of $\C^m$ 
with its defining equations. We introduce two new variables $M,L$ and extra equations which identify these two variables as the eigenvalues of $\rho(m)$ and $\rho(\ell)$ respectively (write $m, \ell$ as words in the generators and look at the upper-left corner). Now that we have augmented $R_T(\Gamma)$ to some algebraic 
set in $\C^{m+2}$, consider the projection to the $M,L$ coordinates into $\C^2$. The closure of the image of this projection is defined by polynomials which can be obtained 
via considering resultants and eliminations from the defining equations of the augmented variety. 

It is known that the A-polynomial of the figure 8-knot complement is given by
$$-M^4+L(1-M^2-2M^4-M^6+M^8) - L^2M^4.$$  
See \cite{murakami} for an explicit derivation.

\subsection{A Topological Application of the A-polynomial}
Let $S$ be a $\pi_1$-injective surface with boundary in a compact 3-manifold $N$ with torus boundary. Fixing a basis $m,\ell$ for $\pi_1 \partial N$, the boundary 
$\partial S$ is equal to $p m + q \ell$ for coprime integers $p,q$ with $p/q \in \Q^*$. The {\em slope of $S$} is the rational number $p/q$. In \cite{CCGLS}, it is shown that 
the slopes of the sides of the Newton polygon of the A-polynomial of $N$ (with respect to the chosen basis $m,\ell$) are among the slopes of incompressible surfaces with boundary in $N$. 

Recall that if $P(x,y) = \sum c_{i,j} x^iy^j$ is an integral polynomial, then the {\em Newton polygon} of $P$ is the convex hull in $\R^2$ of the set 
$\{(i,j) | c_{i,j} \not= 0 \}$. 

\subsection{Gluing Variety and the H-Polynomial}
In this section we will consider the so-called ``gluing variety" which leads to the H-polynomial\footnote{We warn the reader that this is not the standard terminology.} which is for most purposes the same as the A-polynomial. Strictly speaking, the H-polynomial gives the $\pslc$ version of the A-polynomial, see \cite{champanerkar}. The advantage of the H-polynomial is that it is easier to compute compared to the A-polynomial.
  
Let $N$ be a compact 3-manifold $N$ with torus boundary with a fixed basis $\langle m,\ell \rangle \simeq \pi_1 \partial N$. Assume that we are given an ``ideal triangulation", of $N$, that is, we are given tetrahedra $\Delta_1, \hdots, \Delta_n$ such that 
\begin{equation} \label{triangulation}
int(N) := N \setminus \partial N = \bigcup_{j=1}^n \Delta_j^\circ
\end{equation}
where the $\Delta_j^\circ$ denote $\Delta_j$ with vertices taken out. For example, Thurston proved that the complement of the Figure 8-knot admits an ideal triangulation with two tetrahedra.  

\begin{figure}[htbp]
\begin{center}
\includegraphics[width=0.75\textwidth]{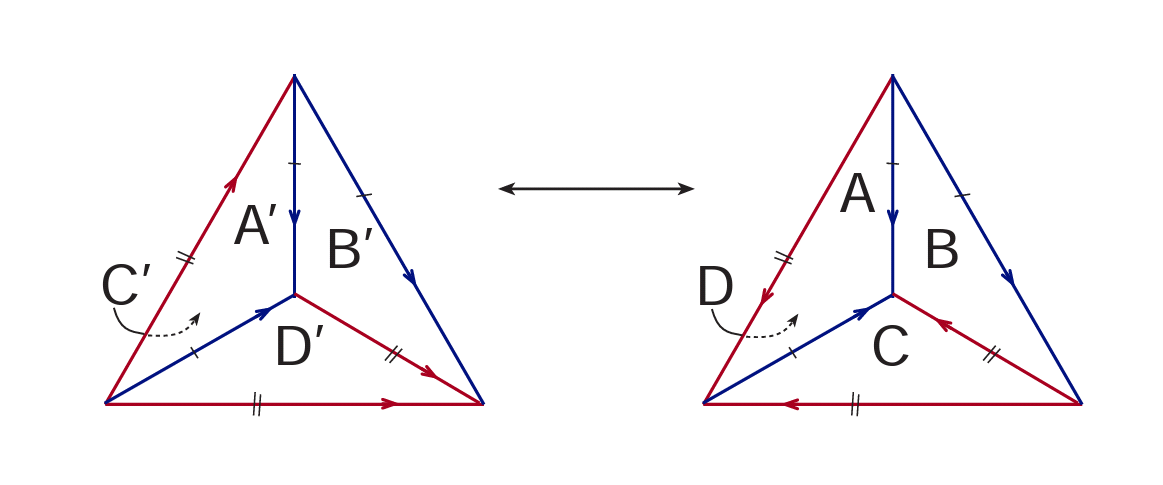}
\end{center}
\caption{Ideal triangulation of the Figure 8-knot complement by two tetrahedra. The image is taken from \cite{hodgson-slides}.}\label{figure8}
\end{figure}

Following Thurston, we would like to investigate the possible hyperbolic structures on $int(N)$ by trying to realize the above ideal triangulation inside the hyperbolic 3-space $\h$.  

Let $\Delta$ be an ideal tetrahedron in $\h$, that is, a tetrahedron whose four vertices lie on $\partial \h = \pc$. Since $\slc$ acts 3-transitively on $\pc$, we can and do apply an isometry to $\Delta$ so that its new vertices are $0,1,\infty,z$ for some $z \in \pc \setminus \{ 0, 1 ,\infty \}$. We will call this $z$ the ``shape parameter" and will denote the associated tetrahedron $\Delta(z)$. The isometry given by the matrix $\left ( \begin{smallmatrix} -1 & 1 \\ -1 & 0 \end{smallmatrix} \right )$ acts as a cyclic permutation on the set $\{ 0 ,1 ,\infty \}$. Thus the tetrahedra $\Delta(z)$, $\Delta(1-1/z)$ and $\Delta(1/(1-z))$ are isometric. Up to replacing by $1-1/z$ or $1/(1-z)$, the shape parameter $z$ of 
$\Delta$ is unique. 

We would like to find $(z_1, \hdots ,z_n)$ such that if we embed $\Delta_1, \hdots , \Delta_n$ in $\h$ as ideal tetrahedra with shape parameters $z_1, \hdots ,z_n$, we will be able to glue them in the way dictated by (\ref{triangulation}) inside $\h$. It turns out that in order to be able to do this, the tuple $(z_1, \hdots ,z_n)$ 
needs to satisfy the $n$ equations, called the {\bf gluing equations}, of the following form
\begin{equation} \label{gluing}
\prod_{i:=1}^n z_i^{a_{i,j}}(1-z_i)^{b_{i,j}} = \epsilon(j), \ \ \ \ \ \ 1 \leq j \leq n
\end{equation}
where $\epsilon(j)= \pm1$ and $a_{i,j}, b_{i,j}$ are integers. In essence, these equations arise from the fact that around every edge, the gluing should be such that 
the angles should sum up to $2\pi$. It turns out that one in fact only needs $n-1$ of these $n$ equations, see \cite{NZ, champanerkar}.

Following \cite{champanerkar}, we define the {\bf gluing variety} of $N$, with respect to the fixed triangulation, as the affine algebraic set given by
 $$G(N) = \{ (z_1, \hdots, z_n,t) \in \C^{n+1} \mid (z_1, \hdots, z_n) \ \textrm{satisfies (\ref{gluing}) and} \ t\prod_{i:=1}^n z_i(1-z_i)=1  \}.$$
The variable $t$ and the equation involving $t$ are there only to ensure that the $z_i \not= 0,1$. It turns out that there is regular map from the gluing variety 
into the $\pslc$-character variety, see \cite{champanerkar, dunfield-2}.

If $(z_1, \hdots ,z_n)$ solves the gluing equations, then we can realize the triangulation (\ref{triangulation}) of $\int(N)$ inside $\h$ using 
$\Delta(z_1), \hdots, \Delta(z_n)$. This puts a hyperbolic metric on $int(N)$. It turns out that the completeness of the resulting metric can analyzed via 
a consideration of the cusp and in particular the metric is complete if and only if $(z_1, \hdots ,z_n)$ solves two more equations of the form 
\begin{equation} \label{complete}
\prod_{i:=1}^n z_i^{c_i}(1-z_i)^{d_i} = 1, \ \ \ \ \ \ \prod_{i:=1}^n z_i^{e_i}(1-z_i)^{f_i} = 1.
\end{equation}
These two equations determine the squares of eigenvalues of the meridian $m$ and longitude $\ell$ in the holonomy representation of $\pi_1N$ to $\pslc$ as rational functions in 
$z_i$'s. 

For example, for the Figure 8-knot complement with the two tetrahedra triangulation mentioned above, the gluing variety is determined by the equation
$z_1(1-z_1)z_2(1-z_2)=1$. This gives an elliptic curve over $\Q$ with conductor 15. Moreover, the completeness equations are $z_1(1-z_1)=1$ and $z_2(1-z_1)=1$. There is essentially a unique solution to these three equations given by $z_1=z_2=\frac{-1+\sqrt{-3}}{2}$, giving the famous complete hyperbolic structure first found by Riley via a consideration of the representation variety.   

Just like we did for the character variety, we would like to obtain a plane curve from the gluing variety by considering the behaviour at the boundary. One can define 
a holonomy map $G(N) \rightarrow \C^* \times \C^*$ taking a point $[z]$ on $G(N)$ to the squares $(x,y)$ of the eigenvalues of $m,\ell$ under the $\pslc$-representation $\rho_{[z]}$ associated to $[z]$. This map gives a curve in $\C^* \times \C^*$ which is usually called the {\em holonomy curve}. 
Its defining polynomial $H(X,Y)$, let's call it the {\bf H-polynomial}, is a factor of the $\pslc$-version of the A-polynomial of $N$, see \cite{champanerkar, dunfield-2} for details. 
 
To compute the $H$-polynomial in practice, we can proceed as we did for the A-polynomial. 
Let us introduce two variables $X,Y$ and consider the following versions of the completeness equations stated in (\ref{complete})
\begin{equation} \label{modified_complete}
\prod_{i:=1}^n z_i^{c_i}(1-z_i)^{d_i} = x, \ \ \ \ \ \ \prod_{i:=1}^n z_i^{e_i}(1-z_i)^{f_i} = y.
\end{equation}
Consider the augmented gluing variety $\overline{G}(N)$ in $\C^{n+3}$ given by the equations that define the gluing variety together with the two equations given in 
(\ref{modified_complete}). The closure of the image of the projection to the $x,y$ coordinates from $\overline{G}(N)$ to $\C^2$ gives us the holonomy curve. 

To illustrate, consider the case of the Figure 8-knot complement $N$. Consider the equation $z_1(1-z_1)z_2(1-z_2)=1$ together with $z_1(1-z_1)=y$ and $z_2(1-z_1)=x$. 
To compute the holonomy curve, denoted $\mathcal{H}$, we eliminate $z_1,z_2$ and find the equation $H(x,y)=0$ where  
$$H(x,y)=y(x^4-x^3-2x^2-x+1)+y^2x^2+x^2.$$ 
Comparing with the A-polynomial, we see that $H(M^2,L)=-A_N(M,L)$.


\section{Mahler Measure}
Given a nonzero Laurent polynomial $P \in \Z[x_1,x_1^{-1} \hdots, x_n,x_n^{-1}]$, the {\bf Mahler measure}\footnote{This is actually called the ``logarithmic Mahler measure", however since we will be 
only working with this, we will simply drop the word ``logarithmic".} of $P$ is defined by
$$m(P) := \int_0^1 \int _0^1 \ln |P(e^{2\pi i \theta_1}, \hdots, e^{2\pi i \theta_n})| d\theta_1 \cdots d\theta_n.$$
Thus the quantity $e^{m(P)}$ gives the geometric mean of $|P|$ on the $n$-torus $S^1 \times \hdots \times S^1$.

When $P \in \Z[x]$, the Mahler measure is intimately related to the notion of ``height" for algebraic numbers. To illustrate, let us start with the so called 
Jensen's Lemma, see \cite[Lemma 1.9]{EW-book} for a proof. 

\begin{lemma}(Jensen's Lemma) \label{Jensen}
For any $\alpha \in \C$, 
$$\int_0^1 \ln |\alpha - e^{2\pi i \theta}| d\theta = \ln^+ |\alpha|,$$
where $\ln^+ \lambda$ denotes $\ln \textrm{max}\{ 1, \lambda \}$.
\end{lemma} 
An application of Jensen's Lemma \ref{Jensen} gives that for $P(x)=a_0 \prod_{j=1}^d (x-\alpha_j)$, we have 
$$m(P) = \ln |a_0| + \sum_{j=1}^d \ln^+ |\alpha_j|.$$

It is a good place to mention the famous question of Lehmer at this point even though we are not interested in it in these lectures. 
\begin{question}(Lehmer) Is $0$ a limit point of the set $\{ m(P) : P \in \Z[x] \}$ ?
\end{question}
The survey \cite{smyth-survey} provides an excellent panaroma of the many aspects of Mahler measure. 

In 1981, Smyth \cite{smyth} proved several elegant formulae one of which is the following:
\begin{equation}\label{smyth}
m(1+x+y) = \dfrac{3\sqrt{3}}{4 \pi} L(\chi_{-3}, 2) = L'(\chi_{-3},-1)
\end{equation}
where 
$$L(\chi_{-3}, s)= \sum_{n \geq 1} \left ( \dfrac{-3}{n} \right ) \dfrac{1}{n^s}$$
is the $L$-function of the character $\chi_{-3}$ of $\Q(\sqrt{-3})$ with conductor $3$. The second equality in formula (\ref{smyth}) follows from the functional equation 
of $L(\chi_{-3},s)$. See the notes of Nuno's talk \cite{nuno} in this mini-workshop for details. The mysterious appearance of an $L$-function in Smyth's formula was explained by work of Deninger in 1997. In \cite{deninger}, Deninger showed
how to interpret $m(P)$ as a {\it Deligne period} of the {\it mixed motive} associated to the variety $V$ given by $P=0$ when $P$ does not vanish on the $n$-torus. Roughly speaking, the mixed motive sits inside the cohomology of $V$ and it follows from conjectures of Beilinson that the determinant of the (Deligne) period matrix of the mixed motive should be related to the value at $s=2$ of the $L$-function associated to $V$. In particular, when $V$ is one-dimensional, one expects, under the Beilinson conjecture, that $m(P)$ 
is directly related to $L(V,0)$. In \cite{deninger}, Deninger came up with the precise prediction that for some $c \in \Q^*$
\begin{equation} \label{deninger}
m(1+x+x^{-1}+y+y^{-1})= c \cdot \dfrac{15}{4\pi^2}L(E,2)= c \cdot L'(E,0)
\end{equation}
where $E$ is the elliptic curve of conductor $15$ given by the projective closure of $1+x+x^{-1}+y+y^{-1} = 0$. Numerical computations carried out by Boyd in \cite{boyd-exp} show that, up to more than 
$50$ digits, the above prediction holds with $c=1$. Recently Rogers and Zudilin proved the equality, see \cite{RZ-15}.
 
We will not go into Deninger's work on which there is a significant amount of interesting work, instead we refer the reader to \cite{deninger, RV-notes, RV-1} for details and 
to \cite{RZ-12, Z-13} and the references in there for some recent progress. We will rather focus on the connections between Mahler measures of A-polynomials and hyperbolic volumes of 3-manifolds, which were discovered by Boyd during the experiments in \cite{boyd-exp}.

\subsection{Mahler Measures of H-Polynomials}

Let $N$ be a compact $3$-manifold with torus boundary. Let $H(x,y)$ denote the H-polynomial of $N$. We have
$$m(H)=\int_0^1 \int _0^1 \ln |H(e^{2\pi i \theta_1},  e^{2\pi i \theta_2})| d\theta_1 d\theta_2$$ 
Let us write $H$ as a polynomial of $x$: $H = a_0(y) \prod_{j=1}^d (x- a_j(y))$
where $a_j(y)$ are rational functions in the variable $L$.  Applying Jensen's formula, we obtain  
$$m(a_0(y)) + \sum_{j=1}^d  \int_0^1  \ln^+ |a_j(y)| d\theta = m(a_0(y)) + \dfrac{1}{2\pi} \sum_{j=1}^d   \int_0^{2\pi}  \ln^+ |a_j(t)| dt.$$
If all roots of $a_0(y)$ lie on the unit circle (such a polynomial is called {\em cyclotomic}), then we have $m(a_0(y))=0$ and it follows that 
$$2 \pi \cdot m(H) = \sum_{j=1}^d   \int_0^{2\pi}  \ln^+ |a_j(t)| dt.$$
Let ${\rm Vol} : \mathcal{H} \rightarrow \R$ be the volume function given by $(x,y) \mapsto \sum_{i=1}^n D(z_i)$ where $(z_1, \hdots, z_n)$ is the solution of the gluing equations determined by $(x,y)$. In \cite{hodgson, dunfield-1}, it is proven that 
$$d{\rm Vol} = -2( \ln|x| d(arg \ y) + \ln|y| d(arg \ x) ).$$
Observe that if $x=e^{i t}$ is on the unit circle, then $d Vol= \ln|y|dt$. 
It follows now that, see \cite{boyd-millennium, BRV-dilog2}, 
$$2 \pi \cdot m(H) = \int_\gamma d{\rm Vol}$$
where $\gamma$ is an oriented path (with possibly many disconnected components) in the intersection of $\mathcal{H}$ with 
$\{ (x,y) \in \C^* \times \C^* \mid |y| = 1, |x| \geq 1 \}$. Now using Stokes Theorem, we get
\begin{equation} \label{formula}
2 \pi \cdot m(H) = \int_{\partial \gamma} {\rm Vol} = \sum_j {\rm Vol}((x_j,y_j)^1) - {\rm Vol}((x_j,y_j)^2)
\end{equation}
where $(x_j,y_j)^1,(x_j,y_j)^2 \in \mathcal{H}$ denote the boundary points of components of $\gamma$. Thus the Mahler measure is given as a sum 
of volumes of $N$ under several different hyperbolic metrics, possibly including the complete one. This gives a conceptual explanation for the 
many numerical examples that Boyd found implying a connection between the Mahler measure of the $A$-polynomial and the volume of the relevant hyperbolic 3-manifold. 

Let us close with our running example: the Figure 8-knot complement. It turns out that in this case (\ref{formula}) involves only the complete hyperbolic metric on $N$ and 
we have $$\pi \cdot m(H) = {\rm Vol}(N).$$


\end{document}